\documentstyle[12pt]{amsart}
\def\F{\hfill $\Box$}
\def\@{$@$}

\def\gz{\ifmmode{Z\hskip -4.8pt Z}
    \else{\hbox{$Z\hskip -4.8pt Z$}}\fi}
\newtheorem{thm}{Theorem}[section]
\newtheorem{prop}[thm]{Proposition}
\newtheorem{lemma}[thm]{Lemma}

\newtheorem{defn}[thm]{Definition}
\newtheorem{cor}[thm]{Corollary}
\newtheorem{exam}[thm]{Example}

\title{SOME DEFINITION OF THE ARTIN EXPONENT OF FINITE GROUPS}
\date{}
\author{K. K. Nwabueze}
\thanks{Research at MSRI is supported by USA National Science Foundation grant 
DMS--9022140}
\address{\hskip -\parindent MSRI and
Department of Maths and Computer Science,
University of Antwerpen (UIA),
2610 Wilrijk, Belgium}
\email{kenneth@@ua.uia.ac.be}
\begin{document}

\begin{abstract}
The Artin exponent induced from cyclic subgroups of finite groups was studied
extensively by T.Y. Lam in \cite{La}.  A Burnside ring theoretic version of
the results in \cite{La} for $p$-groups was given in \cite{Nw}. Here we shall
be interested in looking at the Artin exponent induced from the elementary 
abelian subgroups of finite $p$-groups using some results of A. Dress in \cite{Dr}. 
\end{abstract}    

\maketitle

\section{Introduction}
The Artin exponent induced from cyclic subgroups of finite groups was studied
extensively by T.Y. Lam in \cite{La}. A new setting for the description of 
these exponents in the sense of \cite{La} using the Burnside ring theoretic 
methods was given in \cite{Nw}. In this paper we are interested in defining a 
similar concept for elementary abelian subgroups of finite groups. 
More precisely, using basically some results of A. Dress concerning the 
Burnside rings (see\cite{Dr1},\ \cite{Dr2},\ \cite{Dr}), we shall prove the 
following results (see chapter three for the explanation of the notation).
\medskip

{\bf Main result}
Let $G$ be a finite $p$-group and $e(G)$ the Artin exponent induced from the 
elementary abelian subgroups of $G.$ One has the following:\\
(a)If $G$ is abelian, then $e(G) = |G : {\overline U}|,$ 
where ${\overline U}$ is the maximal elementary
abelian subgroup of $G.$\\
(b) If $G$ is a quaternion or dihedral group then $e(G) = 2$ and $e(G) = 4$ if
$G$ is the semidihedral group.\\
(c) In all other cases  $e(G) = |G|/p. $ 
\medskip

The above results seem to suggest that the invariant $e(G)$ gives an 
interesting numerical measure of the deviation of $G$ from being an elementary
abelian group. Let me mention where this study fits in. First, the method we 
are going to adopt here is Burnside ring theoretic and hence affirms the 
utility of the Burnside ring in the 
representation theory of finite groups. 
Second, as is known, the idea of the Artin exponent appeals to number theory
-- which in themselves have aesthetic values. Finally, exponents of this type 
have been known to be very useful in some formulation of some induction 
theorems for various functors on the category of finite groups -- see for 
instance \cite{Sw}.
\medskip

This paper is arranged as follows.
In section 2, we collect some well known definitions and basic results about
the Burnside rings. In section 3, we form our definition and provide some
particular examples. In section 4, we prove some useful results from which 
the main result finally follows in section 5.

\medskip

\section{The Rings $\Omega (G)$ and $\tilde{\Omega}(G)$}
Let $G$ be a finite group. A $G$-set is a finite set on which $G$ acts from the
left by permutation. The set of isomorphism classes of finite left $G$-sets
form a commutative semi-ring $\Omega^{+}(G)$ with addition induced by disjoint 
union and multiplication induced by cartesian product with diagonal actions. 
The Grothendieck ring of $\Omega^{+}(G)$ is denoted by $\Omega(G)$ and called
the {\it Burnside ring} of $G.$ 
Let $G^{sub/\sim}$ denote the set of conjugate classes of subgroups of $G.$  
For any $G$-set $X$ let $[X]$ be the image in $\Omega (G)$ and write $(U)$ for
the set of subgroups conjugate to $U$ in $G.$ 
Then $\Omega(G)$ is additively a free abelian group with basis $[G/U]$ where 
$(U)$ is taken over the set $G^{sub/\sim}.$ 
For any $G$-set $X$ and a subgroup $U$ of $G,$ let
$X^{U} := \{x \in X\ |\ ux = x\ \forall u \in U \},$ 
be the set of 
$U$-invariant elements of $X.$ Since the map $|X| \longrightarrow |X^{U}|$ 
preserves sums and products, it extends to a ring homomorphism 
$$\chi_{U}:\ \Omega(G) \longrightarrow Z : [X] \longrightarrow |X^{U}|$$ 
from $\Omega(G)$ to $Z$ (the integers). 
Now for every $(U) \in G^{sub/\sim}$ let $Z_{(U)}$ be the ring
isomorphic to $Z.$ We define the {\it ghost ring} $\tilde{\Omega}(G)$ as
follows: 
$$\tilde{\Omega}(G)\ :=\ \mathop\Pi\limits_{(U) \in G^{sub/\sim}} Z_{(U)}.$$
So $\tilde{\Omega}(G)$ is a ring by pointwise multiplication.
It is well known (see \cite{Dr}) that the product map 
$$\chi := \mathop\Pi\limits_{(U) \in G^{sub/\sim}}  \chi_{(U)}
:\ \Omega(G) \longrightarrow \tilde{\Omega}(G)$$ is injective. One also has the
following important result.

\begin{thm} 
Identifying $ \Omega(G) $ with its image in $ \tilde{\Omega}(G) $ with respect 
to the map $ \mathop\Pi\limits_{(U) \in G^{sub/\sim}} \chi_{(U)}$
 one has that $ \{ n \in Z\ |\  n \cdot \tilde{\Omega}(G) \subseteq \Omega(G) \}\ =\ |G|\cdot Z.$
\end{thm}
\medskip

{\bf Proof:} see \cite{Dr}

\begin{lemma} For every $x \in \Omega(G)$ one has 
$\mathop\sum\limits_{g \in G}\chi_{<g>}(x)\   \equiv\  0 (\mbox{mod}\ |G|),$ where $<g>$ denotes the
cyclic group generated by $ g. $ (This relation is called  the Cauchy-Frobenius-Burnside relation).
\end{lemma}
\medskip

{\bf Proof:} see \cite{DSY}
\bigskip

Because for any two subgroups $U$ and $V$ in  $G$ one has that 
$\chi_{U} = \chi_{V}$ if and only if $U \mathop{\sim}\limits^{G} V$ ($U$ and $V$ are conjugate in $G$) and for $x$ and  $y$ in $\Omega(G)$ one has 
$\chi_{U}(x) = \chi_{U}(y)$ for all subgroups $U$ in $G$ 
if and only if $ x = y $ it follows that we can identify each 
$x \in \Omega(G)$ 
with the associated map $ U\ \longrightarrow \chi_{U}(x),$ from 
the set of all subgroups of $G$ into $Z.$ We shall denote this associated 
map also by $x.$ So we can consider $\Omega(G)$ in a 
canonical way as a subring of its ghost ring $ \tilde{\Omega}(G).$

\begin{thm} 
Let $x \in \tilde{\Omega}(G).$ Then $x \in \Omega(G)$ if and only 
if for every $U\ \unlhd\  V\ \leq\  G$ with $(V\ :\ U)$ a power of a prime
one has 
$$ \mathop\sum\limits_{vU \in V/U} x(<v\ ,\ U>)\ \equiv\ 0( mod\ (V : U)).$$
\end{thm}
\medskip

{\bf Proof:} see \cite{Dr2}
\bigskip

For more details on the Burnside ring see \cite{Dr}, \cite{Di}.
\medskip

\section{The Artin Exponent}

We now give the definition of the Artin exponent in terms of the
Burnside ring. This arises in the following context.

\begin{defn}
Let ${\cal U}$ denote the family of all elementary abelian subgroups of $G.$ 
Let $b_{\cal U}$ be the element of ${\tilde \Omega}(G)$ defined as follows. If
$(U)$ is the conjugacy class of the subgroup $U$ of $G$ then
$$ b_{\cal U}(U)\  :=
\begin{cases}
1& \text{ if  $ U \in {\cal U}$}\\
0& \text{ if $U \not\in {\cal U}$}
\end{cases} .$$ 

We call the integer 
$e(G)\ :=\ min(n \in N\ |\ n \cdot b_{\cal U} \in \Omega(G))$ 
the Artin exponent of $G.$
\end{defn}

\begin{cor}
$|G| \cdot b_{\cal U}(U) \in \Omega(G).$ 
\end{cor}

{\bf Proof:} Follows from Theorem 2.1. \F
\medskip

\begin{cor} 
$e(G)$ divides the order of $G.$
\end{cor}

{\bf Proof:} Put $e := e(G)$ and $b := b_{\cal U}(U).$ 
Write $|G| = qe + r$ with $0 \leq r < e.$
Now\\ $r \cdot b\ =\ (|G|- qe) \cdot b\ =\ |G| \cdot b - (qe) \cdot b\ \in\ 
\Omega(G).$ Hence $r = 0$ in view of the minimality of $e.$     \F
\bigskip

To motivate the general algebraic procedure to be followed in the coming
sections we take some particular examples.

\begin{exam}
For a finite group $G,$ let ${\cal S}(G)$ denote the set of all subgroups
of $G.$ Now let $G$ be a cyclic $p$-group of order $p^{\gamma}.$ That is\\ 
$G = <g>,\ |G| = p^{\gamma}$ and ${\cal S}(G) = \{U_{0},\ U_{1},\ U_{2},\ .\ .\ , U_{\gamma-1}\},$
where\\ 
$U_{0} = <g^{p^{\gamma}}> ;\ U_{1} = <g^{p^{\gamma-1}}> ;\ .\ .\ .\ .\ ;\ U_{\gamma-i}\ =\ <g^{p^{i}}>\ \mbox { and }\\ |U_{j}| = p^{j}.$ The family ${\cal{U}}$ of elementary abelian subgroups of $G$ consists
of the subgroups\ $U_{0}$\ and\ $U_{1}$. From definition 3.1, we have that\\
$e(G) \cdot b_{\cal U} = e(G) \cdot (1,\ 1,\ 0,\ .\ .\ .\ ,0)\ =\ (e(G),\ e(G),\ 0,\ .\ .\ .\ ,0) \in \Omega(G).$  To finish this example we need the following 
lemma.

\begin{lemma} Let $G$ denote a cyclic $p$-group of order $p^{n}.$ For every
$i \in \{0,\ .\ .\ .\ ,n\}$  let $U_{i}$ denote the unique subgroup of order 
$p^{i}$ of $G.$  Then given 
$x(U_{i}) \in \tilde{\Omega}(G)$ one has that $x(U_{i}) \in \Omega(G)$ 
if and only if 
$$p^{i} \cdot x(U_{i})\ +\ \mathop \sum \limits_{j=i+1}^{n}(p^{j}-p^{j-1}) \cdot x(U_{j})\ \ \equiv\ 0(mod\ p^{n}), $$ where $ |U_{j}|\ =\ p^{j}.$
\end{lemma}
\medskip

{\bf Proof:} For $ x\ \in\ \Omega(G)$ and $U\ \unlhd\ V\ \leq\ G$ we have that 
$$\sum_{vU \in V/U} x(<v\ ,\ U>)\ \equiv\ 0\mbox{ (mod }V : U).$$ 
Since  $G$ is abelian, $U_{i}\ \unlhd\ G\ \forall i \in \{0,\ .\ .\ .\ ,n\}.$ 
Hence by Theorem 2.3 
$$ \sum_{gU_{i} \in G/U_{i}} x_{i}(<g\ ,\ U_{i}>)\ \ \equiv\ \ 0(mod\ p^{n-i}).$$ 
Since $x(<g\ ,\ U_{i}>)\ =\ x(U_{j})$ if and only if $<g\ ,\ U_{i}>\ =\ U_{j}$
we have\\ 
$ \mathop\sum\limits_{gU_{i} \in G/U_{i}} x(<g\ ,\ U_{i}>)\ \ =\ \ \mathop\sum \limits_{j=i}^{n} x(U_{j}) \cdot \# \{ gU_{i}\ \in\ G/U_{i}\ \ |\ <g\ ,\ U_{i}>\ =\ U_{j}\ \}.$ 
Since we have that\ $ <g\ ,\ U_{i}>\ \ =\ \ <g>\ \ $\ if\ $g\ \not\in\  U_{i}$
then
\begin{eqnarray*}
& &\ \ \ \ \ \ \ \ \ \ \mathop\sum\limits_{j=i}^{n} x(U)_{j}\ \cdot \# \{ gU_{i}\ \in\ G/U_{i}\ \ |\  \ <g\ ,\ U_{i}>\ =\ U_{j} \}\\ 
& &\ \ \ =\ \ \ \ x(U)_{i} + \sum_{j=i+1}^{n} x(U_{j}) \cdot \frac{p^{j}-p^{j-1}}{p^{i}}\ \ \ \equiv\ \ \ 0(\mbox{mod }p^{n-i})\\
& &\Longleftrightarrow\ \ \ \ x(U_{i}) \cdot p^{i}\ +\ \sum_{j=i+1}^{n} x(U_{j}) \cdot (p^{j}-p^{j-1})\ \ \ \equiv\ \ \ 0( \mbox{mod }p^{n}).
\end{eqnarray*}                                             
\F                                                                 
\bigskip

Now we finish our example by observing from the above lemma that for           $i=1$ we have  $e(G) \cdot p\ \ \equiv\ \ 0(mod\ p^{\gamma})$ and for 
$i = 0$ we have 
$$e(G)\ +\ e(G) \cdot (p - 1)\ \ \equiv\ \ 0(mod\ p^{\gamma})\ \ \Longleftrightarrow\ \ \  p^{\gamma-1}\ \ |\ \ e(G).$$  This implies $e(G)\ =\ p^{\gamma-1}.$  \end{exam}                        
\medskip

\begin{exam}
Consider ${\cal U}$ a family of elementary abelian subgroups of a finite group
$G.$ Then, from lemma 3.5, one has that $e(G) = 1$ if and only if 
${\cal U} = {\cal S}(G).$ In particular, if $G$ is elementary abelian then $e(G) = 1.$
\end{exam}
\medskip

\begin{exam}
Let $G$ be the quaternion group of order $8.$ Then it is to see by direct computation that $e(G) = 2.$ 
\end{exam}
\bigskip

We now generalize the above examples in the following section.

\section{Results}
We start this section with the following observation.

\begin{lemma}  For ${\cal{U}}$ a family of elementary abelian 
subgroups of a finite group $G$ let  $e_{\cal U}(G)$ be the Artin exponent 
with respect to the family ${\cal U}.$  Assume that $ e_{\cal U}(G) = 1,$ 
then for all $U , V \leq G$ with $U \unlhd V$ and $(V : U)$ a prime power, one 
has that $U \in {\cal{U}}$ if and only if $V \in {\cal{U}}.$
\end{lemma}
\medskip

{\bf Proof:}  Assume $(V : U)\ =\ p^{\alpha}$ for some $\alpha \geq 0.$ 
Then $x  \in \Omega(G)$ implies 
$x(U) \equiv x(V)(mod\  p).$ Therefore $ U \in {\cal{U}}$ implies 
$x(V) \equiv  1(mod\ p),$ hence $x(V) \ne 0.$  We have $x(V) = 1$ and so 
$V \in {\cal{U}}.$ Conversely $U \not\in {\cal{U}}$ implies 
$x(V)\ \equiv\ 0(mod\ p), $ hence\ $x(V)\ \ne\ 1$ which implies  
$x(V)= 0,$ that is $V \not\in {\cal{U}}.$                                    \F

\begin{lemma} 
If $G$ is solvable, then $x \in \Omega(G)$ implies that 
${\cal U} = {\cal S}(G)$ or ${\cal{U}} = \emptyset$ (the empty set).
\end{lemma}
\medskip

{\bf Proof:} If $G$ is solvable then for every $U$ in $G$ there exists a 
sequence of subgroups $U\ =\ U_{0},\ U_{1},\ .\ .\ .\ .\ , U_{k} = 1$ such 
that $U_{i}$ is normal in\ $ U_{i-1} $\ for\ $ i = 1,\ .\ .\ .\ .\ .\ ,k $\ and\ $ (U_ {i-1} : U_{i}) $\ is a prime power or $1.$ Hence if $ {\cal{U} }\ \ne\  \emptyset,$ then by (4.1) one has that 
$1 \in {\cal{U}} $ and $ U \in {\cal{U}} $ for all $ U \leq G. $   
\F
\medskip

\begin{cor} If $G$ is  a finite $p$-group, then $e(G) = 1$ if and only if 
${\cal{U}} = {\cal S}(G).$ 
\end{cor}
\medskip

{\bf Proof:} Follows from 4.2. 
\F
\medskip

\begin{prop} 
Let $G$ be a finite noncyclic abelian $p$-group. Then one has 
$e(G) = (G : \overline{U}),$
where $\overline{U} = \{ g \in G\ \ |\ \ g^{p}\ =\ 1 \}$ 
is the maximal elementary abelian subgroup of $G.$ 
\end{prop}
\medskip

{\bf Proof:} Let ${\cal U}$ be the family of all elementary abelian subgroups
of $G$ and let $U \in {\cal U}.$ Let $e(G)$ and $b_{\cal U}$ be defined as in
3.1. By definition of $e(G)$ we have $e(G) \cdot b_{\cal U} \in \Omega(G),$ 
hence by Theorem 2.3, 
$$ \mathop\sum\limits_{gU \in G/U} e(G) \cdot b_{\cal U}(<g\ ,\ U>)\ \equiv\ 
0(mod\ |G\ :\ U|).$$
Since $<g\ ,\ U>$ is elementary abelian if and only if $g \in {\overline U},$ 
we have 
\begin{eqnarray*}
& &\mathop\sum\limits_{gU \in G/U} b_{\cal U}(<g\ ,\ U>)\\ & &\ \ =\ |\{ gU \in G/U\ \ |\ \ <g\ ,\ U> \mbox{ is elementary abelian } \}|\\
& &\ \ =\ |\{ gU \in G/U\ \ |\ \ <g\ ,\ U>\ \leq\ {\overline U}\}|\\
& &\ \ =\ |{\overline U}\ :\ U|.
\end{eqnarray*}

Hence it follows that $e(G)$ is the minimal positive integer such that 
$$ e(G) \cdot |{\overline U}\ :\ U|\ \equiv\ 0(mod\ |G\ :\ U|),$$
that is $e(G)\ =\ (G : \overline{U}).$
\F

\medskip

Before we state the next result we need the following (see \cite{Za}).
\medskip

\begin{lemma} Let $G$ be a non-abelian $p$-group of order $p^{n} (p = 2).$ 
If there exists a $ g \in G $  of order $ p^{n-1} $ then  $G$ has one of the 
following  presentations for some $h$ in $G.$\\
(A)  $g^{2^{n-1}} = 1,$\ $h^{2} = g^{2^{n-2}}$,\ $hgh^{-1} = g^{-1},$
where $p=2$ and $n \geq 3,$ -- the quaternion group.\\
(B)  $g^{2^{n-1}} = 1$, $h^{2} = 1$, $hgh^{-1} = g^{-1},$ where $p=2$ and 
$n \geq 3$ -- the dihedral group.\\
(C) $g^{2^{n-1}} = 1$,\ $h^{2} = 1$,\ $hgh^{-1} = g^{1+2^{n-2}},$ where 
$p = 2$ and $n \geq 4.$\\
(D) $g^{2^{n-1}} = 1$,\ $h^{2} = 1 $,\ $hgh^{-1} = g^{-1+2^{n-2}}.$ where
$p = 2$ and $n \geq 4,$ -- the semi-dihedral group.
\end{lemma}
\medskip

\begin{prop} 
Let  $G$ be a finite nonabelian and noncyclic $p$ group. Then one has that\\ 
$e(G) = |G|/p$ unless $G$ is 
the quaternion, the dihedral or the semi-dihedral group in which cases 
$e(G) = 2$ for the quaternion or the dihedral groups and $e(G) = 4$ for the 
semi-dehidral.
\end{prop}
\medskip

{\bf Proof:} As usual let ${\cal U}$ be the family of all elementary abelian 
subgroups of $G$ and let $U \in {\cal U}.$ Let $e(G)$ and $b_{\cal U}$ be 
defined as before. So $e(G) \cdot b_{\cal U} \in \Omega(G)$ and 
$\mathop\sum\limits_{gU \in G/U} e(G) \cdot b_{\cal U}(<g\ ,\ U>)\ \equiv\ 
0(mod\ |G\ :\ U|).$ To find $e(G),$ it suffices to set $|U| = p.$ Now
\begin{eqnarray*}
& &\mathop\sum\limits_{gU \in G/U} b_{\cal U}(<g\ ,\ U>)\\ 
& &\ \ =\ \ |\{ gU \in G/U\ \ |\ \ <g\ ,\ U> \mbox{ is elementary abelian } \}|\\
& &\ \ =\ \ 1\ +\ |\{ gU \in G/U\ \ |\ \ <g\ ,\ U>\ \ \cong\ \ C_{p} \times C_{p} \}|,
\end{eqnarray*}
where $C_{p}$ denote a cyclic group of order $p.$ Let 
$$ {\cal N}(g\ ,\ U)\ :=\ |\{ gU \in G/U\ \ |\ \ <g\ ,\ U>\ \cong\ C_{p} \times C_{p} \}|.$$ 
For $G$ a noncyclic $p$-group, it is easy to see that 
${\cal N}(g\ ,\ U)\ \equiv\ 0(mod\ p)$ if $G$ is not a quaternion, dihedral
or semidihedral group. This implies that $e(G) = |G|/p$ is the
desired minimal number such that $e(G) \cdot (1 + {\cal N}(g\ ,\ U)) \equiv 0(mod\ (G\ :\ U)),$ that is $e(G) = |G|/p.$   
On the other hand if ${\cal N}(g\ ,\ U)\ \not\equiv\ 0(mod\ p)$ then $G$ corresponds to the
nonabelian group of order $8,$ i.e. the quaternion and the dihedral groups and
by easy computation we find that in these cases $e(G) = 2.$  
The semi-dihedral group of order $8$ is just the noncyclic abelian group of
order $8.$ So for the semi-dihedral $e(G) = 4.$  
\F
\medskip

\section{Proof of main result}
Finally we have;

\begin{thm}
Let $G$ be a finite $p$-group. We have the following.\\
(a)If $G$ is abelian, then the Artin exponent $e(G) = |G : U|,$ where $U$ is the maximal elementary abelian subgroup of $G.$\\
(b) If $G$ is a quaternion or dihedral group, then $e(G) = 2,$ and $e(G) = 4$ 
if $G$ is a semidihedral group.\\
(c) In all other cases  $e(G) = |G|/p $ 
\end{thm}

{\bf Proof:} Results now follows from proposition 4.4 and proposition 4.6.

\medskip
{\samepage
\pagestyle{plain}
{}
}

%
%
%
%
%
%
%
%

\allowbreak

\end{document}